\documentclass[12pt]{amsart}
\usepackage {amssymb, adjustbox, enumerate, amsbsy, stmaryrd}
\usepackage {amsmath}
\usepackage[mathscr]{eucal}
\usepackage{amsthm}
\usepackage{graphicx}
\usepackage {amscd}
\usepackage {epic}
\usepackage[alphabetic]{amsrefs}
\usepackage{color}
\usepackage{tikz-cd}
\usepackage[all,color]{xy}
\usepackage{pb-diagram}
\usepackage[mathscr]{eucal}
\usepackage{hyperref}
\hypersetup{
    colorlinks=true, 
    linktoc=all,     
    linkcolor=blue,  
}
\usepackage{enumitem,titletoc}

\usepackage{xcolor}

\newcommand{\black}{\color{black}}

\usepackage[normalem]{ulem}

\DeclareMathAlphabet{\mathpzc}{OT1}{pzc}{m}{it}

\usepackage{amsfonts}
\usepackage{amsmath}
\usepackage{amssymb}

\begin{document}

\newcommand{\ssub}{\subset\joinrel\subset}

\renewcommand{\thefootnote}{\fnsymbol{footnote}}
\renewcommand{\thanks}[1]{\footnote{#1}} 
\newcommand{\starttext}{ \setcounter{footnote}{0}
\renewcommand{\thefootnote}{\arabic{footnote}}}
\renewcommand{\theequation}{\thesection.\arabic{equation}}
\newcommand{\be}{\begin{equation}}
\newcommand{\bea}{\begin{eqnarray}}
\newcommand{\eea}{\end{eqnarray}} \newcommand{\ee}{\end{equation}}
\newcommand{\N}{{\cal N}} \newcommand{\<}{\langle}
\renewcommand{\>}{\rangle}
\def\ba{\begin{eqnarray}}
\def\ea{\end{eqnarray}}
\newcommand{\PSbox}[3]{\mbox{\rule{0in}{#3}\includegraphics{#1}
\hspace{#2}}}

\def\v{\vskip .1in}

\def\al{\alpha}
\def\b{\beta}
\def\c{\chi}
\def\d{\delta}
\def\e{\epsilon}
\def\g{\gamma}
\def\l{\lambda}
\def\m{\mu}
\def\n{\nu}
\def\o{\omega}
\def\f{\varphi}
\def\r{\rho}
\def\si{\sigma}
\def\t{\theta}
\def\z{\zeta}
\def\k{\kappa}

\def\G{\Gamma}
\def\D{\Delta}
\def\O{\Omega}
\def\T{\Theta}

\newcommand{\ddbar}{\sqrt{-1}\partial\dbar}

\def\cA{{\mathcal A}}
\def\cB{{\mathcal B}}
\def\cC{{\mathcal C}}
\def\cD{{\mathcal D}}
\def\cE{{\mathcal E}}
\def\cF{{\mathcal F}}
\def\cG{{\mathcal G}}
\def\cH{{\mathcal{H}}}
\def\cI{{\mathcal I}}
\def\cK{{\mathcal K}}
\def\cL{{\mathcal L}}
\def\cM{{\mathcal M}}
\def\cN{{\mathcal N}}
\def\cO{{\mathcal O}}
\def\cP{{\mathcal P}}
\def\cp{{\mathcal p}}
\def\cR{{\mathcal R}}
\def\cS{{\mathcal S}}
\def\cT{{\mathcal T}}
\def\cV{{\mathcal V}}
\def\cX{{\mathcal X}}
\def\cY{{\mathcal Y}}
\def\A{{\mathbb{A}}}
\def\N{\mathbb N}
\def\Z{{\mathbb Z}}
\def\Q{{\mathbb Q}}
\def\R{{\mathbb R}}
\def\C{{\mathbb C}}
\def\P{{\mathbb P}}
\def\K{{\rm K\"ahler }}

\def\KE{{\rm K\"ahler-Einstein }}
\def\KEE{{\rm K\"ahler-Einstein}}
\def\Rm{{\rm Rm}}
\def\Ric{{\rm Ric}}
\def\Hom{{\rm Hom}}
\def\mod{{\ \rm mod\ }}
\def\Aut{{\rm Aut}}
\def\End{{\rm End}}
\def\osc{{\rm osc\,}}
\def\vol{{\rm vol}}
\def\Vol{{\rm Vol}}
\def\reg{{\rm reg}}
\def\sing{{\rm sing}}
\def\KE{{\rm K\"ahler-Einstein\ }}
\def\PSH{{\rm PSH}}
\def\Ad{{\rm Ad}}
\def\Lie{{\rm Lie}}
\def\ti\tilde
\def\ann{{\rm ann\,}}
\def\u{\underline}

\def\pl{\partial}
\def\na{\nabla}
\def\i{\infty}
\def\I{\int}
\def\p{\prod}
\def\s{\sum}
\def\dd{{\bf d}}
\def\ddb{\partial\bar\partial}
\def\sub{\subseteq}
\def\ra{\rightarrow}
\def\hra{\hookrightarrow}
\def\Lra{\Longrightarrow}
\def\lra{\longrightarrow}
\def\LA{\langle}
\def\RA{\rangle}
\def\L{\Lambda}
\def\diam{{\rm diam}}
\def\Diff{{\rm Diff}}
\def\Mod{{\rm Mod}}
\def\Hilb{{\rm Hilb}}
\def\depth{{\rm depth}}
\def\Ass{{\rm Ass}}
\def\us{{\underline s}}
\def\Re{{\rm Re}}
\def\Im{{\rm Im}}
\def\tr{{\rm tr}}
\def\det{{\rm det}}
\def\half{ {1\over 2}}
\def\third{{1 \over 3}}
\def\ti{\tilde}
\def\un{\underline}
\def\Tr{{\rm Tr}}
\def\Ker{{\rm Ker}}
\def\spec{{\rm Spec}}
\def\supp{{\rm supp}}
\def\Id{{\rm Id}}

\def\pz{\partial _z}
\def\pv{\partial _v}
\def\pw{\partial _w}
\def\w{{\bf w}}
\def\x{{\bf x}}
\def\y{{\bf y}}
\def\tet{\vartheta}
\def\dwplus{\D _+ ^\w}
\def\dxplus{\D _+ ^\x}
\def\dzplus{\D _+ ^\z}
\def\chiz{{\chi _{\bar z} ^+}}
\def\chiw{{\chi _{\bar w} ^+}}
\def\chiu{{\chi _{\bar u} ^+}}
\def\chiv{{\chi _{\bar v} ^+}}
\def\os{\omega ^*}
\def\ps{{p_*}}

\def\hO{\hat\Omega}
\def\ho{\hat\omega}
\def\o{\omega}
\def\KSBA{{\rm KSBA}}
\def\CM{{\rm CM}}
\def\WP{{\rm WP}}

\def\[{{\bf [}}
\def\]{{\bf ]}}
\def\Rd{{\bf R}^d}
\def\Ci{{\bf C}^{\infty}}
\def\pl{\partial}
\def\sq{{{\sqrt{{\scalebox{0.75}[1.0]{\( - 1\)}}}}}\hskip .01in}
\newcommand{\dotcup}{\ensuremath{\mathaccent\cdot\cup}}

\newcommand{\ssubset}{\subset\joinrel\subset}


\newcommand{\mM}{{\mathcal{M}}}
\newcommand{\ov}{\overline}
\newcommand{\tO}{{\widetilde{\Omega}}}
\newcommand{\tU}{{\widetilde{U}}}
\newcommand{\tOt}{{\widetilde{\Omega}_{\cX_t}}}
\newcommand{\tOo}{{\widetilde{\Omega}_{\cX_0}}}
\newcommand{\ks}{{\rm ks}}
\newcommand{\kss}{{\rm kss}}
\newcommand{\kps}{{\rm kps}}
\newcommand{\rdiv}{{\rm div}}
\newcommand{\mU}{{\mathcal{U}}}
\newcommand{\mV}{{\mathcal{V}}}
\newcommand{\tcY}{{\widetilde{\mathcal{Y}}}}
\newcommand{\tcX}{{\widetilde{\mathcal{X}}}}
\newcommand{\cZ}{{\mathcal{Z}}}
\newcommand{\pdJd}{{\frac{\sqrt{-1}}{2\pi}\partial\bar{\partial}}}
\newcommand{\ku}{{\mathfrak{u}}}
\newcommand{\tku}{\widetilde{\ku}}
\newcommand{\kU}{{\mathfrak{U}}}

\newcommand{\Fut}{{\rm DF}}
\newcommand{\gh}{\rm graph }
\newcommand{\MO}{\mathcal{O}}
\newcommand{\bX}{{\bar{\mathcal{X}}}}
\newcommand{\bB}{\mathbf{B}}
\newcommand{\bD}{\mathbf{D}}
\newcommand{\bd}{\mathbf{d}}
\newcommand{\bM}{\mathbf{M}}
\newcommand{\bL}{{\bar{\mathcal{L}}}}
\newcommand{\bT}{\mathbf{T}}
\newcommand{\mY}{{\mathcal{Y}}}
\newcommand{\mL}{{\mathcal{L}}}
\newcommand{\mX}{{\mathcal{X}}}
\newcommand{\bcM}{{\bar{\mathcal{M}}}}
\newcommand{\bcH}{{\bar{\mathcal{H}}}}
\newcommand{\bcY}{{\bar{\mathcal{Y}}}}
\newcommand{\DMR}{{\mathcal{DMR}}}
\newcommand{\DR}{{\mathcal{DR}}}
\newcommand{\mop}{{\mathcal{O}_{\mathbb{P}^1}(1)}}
\newcommand{\mot}{{\mathcal{O}_{\mathbb{P}^1}(-1)}}
\newcommand{\mult}{{\rm mult}}
\newcommand{\mld}{{\rm mld}}
\newcommand{\Supp}{{\rm Supp}}

\newcommand{\AAA}{\mathbb{A}}
\newcommand{\QQ}{\mathbb{Q}}
\newcommand{\GG}{\mathbb{G}}
\newcommand{\PP}{\mathbb{P}}
\newcommand{\RR}{\mathbb{R}}
\newcommand{\CC}{\mathbb{C}}
\newcommand{\ZZ}{\mathbb{Z}}
\newcommand{\kk}{\mathbb{k}}
\newcommand{\sE}{\mathscr{E}}
\newcommand{\sL}{\mathscr{L}}
\newcommand{\sO}{\mathscr{O}}
\newcommand{\sI}{\mathscr{I}}
\newcommand{\sJ}{\mathscr{J}}
\newcommand{\cJ}{\mathcal{J}}
\newcommand{\cU}{\mathcal{U}}
\newcommand{\fg}{\mathfrak{g}}
\newcommand{\fk}{\mathfrak{k}}
\newcommand{\fp}{\mathfrak{p}}
\newcommand{\fsl}{\mathfrak{sl}}
\newcommand{\ddJ}{\sqrt{-1}\partial_J\bar\partial_J}
\newcommand{\id}{\mathrm{id}}
\newcommand{\aut}{\mathrm{Aut}}
\newcommand{\faut}{\mathfrak{aut}}
\newcommand{\dist}{\mathrm{dist}}
\newcommand{\Sing}{\mathrm{sing}}
\newcommand{\chow}{\mathrm{Chow}}
\newcommand{\hilb}{\mathrm{Hilb}}
\newcommand{\Ext}{\mathrm{Ext}}
\newcommand{\Pic}{\mathrm{Pic}}
\newcommand{\DF}{\mathrm{DF}}
\newcommand{\DFn}{\mathrm{DF}^n}
\newcommand{\GH}{\mathrm{GH}}
\newcommand{\CH}{\mathrm{CH}}
\newcommand{\kst}{\mathrm{kst}}
\newcommand{\FS}{\mathrm{FS}}
\newcommand{\SL}{\mathrm{SL}}
\newcommand{\topo}{\mathrm{top}}
\newcommand{\Red}{\textcolor{red}}
\newcommand{\rred}{}
\newcommand{\Blue}{}
\newcommand{\BBlue}{}
\newcommand{\pp}{{\prime\prime}}
\newcommand{\sddb}{{\sqrt{-1}\partial\bar{\partial}}}

\newtheorem{claim}{Claim}[section]
\newtheorem{theorem}{Theorem}[section]
\newtheorem{proposition}{Proposition}[section]
\newtheorem{lemma}{Lemma}[section]
\newtheorem{conjecture}{Conjecture}[section]
\newtheorem{example}{Example}[section]
\newtheorem{definition}{Definition}[section]
\newtheorem{corollary}{Corollary}[section]
\newtheorem{remark}{Remark}[section]

\address{$^*$ Department of Mathematics, Rutgers University, Piscataway, NJ 08854}

\address{$^{**}$ Department of Mathematics and Computer Science, Rutgers University, Newark, NJ 07102}
\address{$^\dagger$Department of Mathematics and Computer Science,  Rutgers University, Newark, NJ 07102}

\centerline{\bf   AN ANALYTIC PROOF OF THE STABLE REDUCTION THEOREM }

\bigskip

\centerline{ \small  JIAN SONG$^*$, JACOB STURM$^{**}$, XIAOWEI WANG$^\dagger$   \footnote{{ Research supported in part by National Science Foundation grant  DMS-1711439, DMS-1609335 and  Simons Foundation Mathematics and Physical Sciences-Collaboration Grants, Award Number: 631318.}}}

 \bigskip
 \medskip

{\noindent \small A{\scriptsize BSTRACT}. \footnotesize $~$~~  The stable reduction theorem says that a family of curves of genus $g\geq 2$ over a punctured curve can be uniquely completed (after possible base change) by inserting certain stable curves at the punctures. We give a new proof of this result for curves defined over $\C$, using the \KE metrics on the fibers to obtain the limiting stable curves at the punctures.   }

\bigskip


\section{Introduction}
 
Let $X_1, X_2, ...$ be a sequence of compact Riemann surfaces of genus $g\geq 2$. 
A consequence of the Deligne-Mumford construction of moduli space is the following. There exists $N>0$ and imbeddings $T_i: X_i\hookrightarrow \P^N$
such that after passing to a subsequence, $T_i(X_i)=W_i\sub\P^N$ converges to a stable algebraic
curve, i.e. a curve $W_\i\sub\P^N$ whose singular locus is either empty or consists of nodes, and whose smooth locus carries a metric of constant negative curvature. The stable reduction theorem \cite{DM} (stated below) is the analogue of this result with 
$\{X_i\,:\, i \in \N\}$ replaced by an algebraic family $\{X_t\,:\, t\in \D^*\}$ where $\D^*\sub\C$ is the punctured unit disk. 

\v

The imbeddings $T_i$ are determined by a canonical (up to a uniformly bounded automorphism) basis of $H^0(X_i, mK_{X_i})$ (here $m\geq 3$ is fixed). We are naturally led to ask: Can one construct the canonical basis defining $T_i$ explicitly? In Theorem \ref{main} we give an affirmative answer to this question.
\v
The main goal of this paper is to give an independent analytic proof of these algebraic compactness results, which is the content of Theorem \ref{1point2}. We start with the Bers compactness theorem, which says that after passing to a subsequence, the $X_i$ converge to a nodal curve in the Cheeger-Colding topology. We then use the technique of Donaldson-Sun \cite{DS} which uses the K\"ahler-Einstein metric to build a bridge between analytic convergence (in Teichmuller space) to algebraic convergence (in projective space). The main difficulty is that unlike the \cite{DS} setting, the diameters of the $X_i$ are unbounded and as a consequence, some of the pluri-canonical sections on $X_\i$ are not members of $L^2(X_\i, \o_{\rm KE})$, so one can't apply the $L^2$-Bergman imbedding/peak section method directly. In order  
to solve this problem, we introduce the ``$\e$-Bergman inner product"
on the vector space $H^0(X_i, mK_{X_i})$, which is defined by the $L^2$ norm on the thick part of the $X_i$ (unlike the standard Bergman inner product which is the $L^2$ norm defined by integration on all of $X_i$)
and we show that for fixed $m\geq 3$  the canonical basis defining $T_i$ is an an orthnormal basis for this new inner product. This establishes Theorem \ref{main} which we then use to prove Theorem 
\ref{1point2} (the stable reduction theorem).
\v 
We start by reviewing the corresponding compactness results for Fano maniolds established by Donaldson-Sun in \cite{DS}.\,
Let $(X_i,\o_i)$ be a  sequence of  K\"ahler-Einstein  manifolds of dimension $n$ with 
$c_1>0$, volume at least $V$ and diameter at most $D$, normalized so
that $\Ric(\o_i)=\o_i$.   The first step in the proof of the Donaldson-Sun theorem is the application of Gromov's compactness theorem which implies that after passing to a subsequence,  $X_i$ converges to a compact metric space $X_\i$ of dimension $n$
in the metric sense, i.e. the Cheeger-Colding (CC) sense. This first step is not not available in the $c_1<0$ case due to the possibility of collapsing and unbounded diameter. Nevertheless, the analogue of this Cheeger-Colding property for Riemann surfaces of genus $g\geq 2$ is available thanks to the compactness theorem of Bers \cite{B}.  

\v
For the second step, Donaldson-Sun construct explicit imbeddings $T_i: X_i\hookrightarrow\P^N$ with the following properties.
 Let $X_i\ra X_\i$
 in the Cheeger-Colding sense as above. 
Then  there is a K-stable algebraic variety $W_\i\sub\P^N$
such that if $W_i=T_i(X_i)$ then 
$W_i\ra W_\i$ in the algebraic sense (i.e. as points in the Hilbert scheme). Moreover,
$ T_\i: X_\i\ra W_\i
$
is a homeomorphism, biholomorphic on the smooth loci, where 

\be\label{Tinf}\hbox{
$T_\i(x_\i)=\lim_{i\to\i} T_i(x_i)$\ \  whenever\ \  $x_i\ra x_\i$.
}
\ee
\v
We summarize this result with the following diagram:

\be\label{diag}
\begin{tikzcd}
  & X_i\ \arrow{r}{T_i}\arrow{d}[swap]{\rm CC} 
  &\  W_i\ \arrow[hookrightarrow]{r}{}\arrow{d}{\rm Hilb} 
  & \ \P^N 
  & {} \\
  & X_\i\ \arrow{r}{T_\i} 
  &\ W_\i\ \arrow[hookrightarrow]{r}[swap]{} 
  & \ \P^N  
  \end{tikzcd}
\ee
 
Here the vertical arrows represent convergence in the metric (Cheeger-Colding) sense and the the algebraic (Hilbert scheme) sense respectively. The  horizontal arrows isomorphisms:  $T_i$ is an algebraic isomorphism, and $T_\i$ is a holomorphic isomorphism. For $1\leq i\leq \i$, the maps $W_i\hookrightarrow\P^N$ are inclusions.
\v
The imbeddings 
$T_i: X_i\ra\P^N$ are the so called
``Bergman imbeddings". This means  $T_i=(s_0,...,s_N)$ where the $s_\al$ form an
orthonormal basis of $H^0(X_i,-mK_{X_i})$ with respect to the Bergman inner product:

\be\label{berg}
\I_{X_i} (s_\al,s_\b)\,\o_i^n\ = \ \d_{\al,\b}
\ee

Here  $m$ is a fixed integer which is independent of $i$ and the pointwise inner product is defined by $(s_\al,s_\b)=s_\al\bar s_\b\o_i^m$.  Since the definition of $T_i$ depends on the choice of orthonormal basis 
$\us=(s_0,...,s_N)$, we shall sometimes write $T_i = T_{i,\us}$ when we want to stress
the dependence on $\us$. 

\v

 Thus we assume $\Ric(\o_i)=-\o_i$ and we wish
to construct imbeddings $T_i: X_i\ra\P^N$ such that the sequence $W_i = T_i(X_i)\sub\P^N$ converges to a
singular  K\"ahler-Einstein variety $W_\i$ with $K_{W_\i}>0$. 
\v
The condition that $W_\i$ is a ``singular  K\"ahler-Einstein variety" can be made precise as follows.
Let $W\sub\P^N$ be a projective variety with $K_W$ ample. The work of Berman-Guenancia \cite{BG} combined with the results of Odaka \cite{0} tell 
us that the following conditions are equivalent.

\begin{enumerate}
\item There is a \K metric $\o$ on $W^{\rm reg}$ such that $\Ric(\o)=-\o$ satisfying the volume condition $\I_{W^{\rm reg}}\o^n=c_1(K_W)^n$.
\vskip .02in
\item $W$ has at worst semi-log canonical singularities.
 \item $W$ is K-stable
\end{enumerate}

We wish to construct $T_i$ in such a way that  $W_\i=\lim_{i\to\i}T_i(X_i)$
has at worst semi-log canonical singularities. In this paper we  restrict
our attention to the case $n=1$. 
\v
Our long-term goal is to generalize the above theorem of \cite{DS} to the case where the $(X_i,\o_i)$ are
smooth canonical models, of dimension $n$, i.e. $X_i$ is smooth and $c_1(X_i)<0$. The proof we present here is designed with that goal in mind. There are other approaches, but this is the one that seems to lend itself most easily to generalization. We have been able to extend the techniques to the case of dimension two, but that will be the subject a future paper.

\v
\begin{remark}\label{fail}
One might guess, in parallel with the Fano setting, that the $T_i: X_i\ra\P^N$ should be
the pluricanonical Bergman imbeddings, that is $T_i=T_{i,\us}$ where 
$\us=(s_0,...,s_N)$ and the $s_\al$
form an orthonormal basis of $H^0(X_i,mK_{X_i})$ with respect to the inner product (\ref{berg}). But as we shall see,
this doesn't produce the correct limit, i.e. $W_\i$, the limiting variety, is not stable. In order to get the right imbedding into projective space, we need to replace $T_{i,\us}$ with $T_{i,\us}^\e$, the so called $\e$-Bergman imbedding, defined below.
\end{remark}
\v
We first need to establish some notation. Fix $g\geq 2$ and $\e>0$. If $X$ is a compact Riemann surface
of genus $g$, or more generally a stable analytic curve (i.e. a Riemann surface with nodes whose universal cover is the Poincar\'e disk) of genus $g$, we define the $\e$-thick part of $X$ to be
$$ X_\e\ := \ \{x\in X\,:\, {\rm inj}_x\geq\e\}
$$
Here ${\rm inj}_x$ is the injectivity radius at $x$ and  the metric $\o$ on $X$  is the unique hyperbolic metric satisfying $\Ric(\o)=-\o$.
It is well known that there exists $\e(g)>0$ such that for all $X$  of genus $g$, and for all $0<\e<\e(g)$, that $X\backslash X_\e$ is a finite disjoint union of holomorphic annuli.
\v
Next we define the ``$\e$-Bergman imbedding" $T_{\us}^\e: X\ra\P^N$. Fix
$0<\e<\e(g)$ and fix $m\geq 3$. For each stable analytic curve of genus $g$,  
 we choose a basis $\us=\{s_0,...,s_{N_m}\}$ of
$H^0(X,mK_X)$
 such that
$$  \I_{X_\e} (s_\al,s_\b)\,\o\ = \ \d_{\al,\b}
$$
Here 
$(s_\al,s_\b)=s_\al\bar s_\b\o_i^{-m}$ is the  usual pointwise inner product.
Such a basis is uniquely determined up to the action of $U(N+1)$.
Let $T_{\us}^\e: X\hookrightarrow\P^{N_m}$ be the map $T_\us^\e=(s_0,...,s_{N_m})$. Let $W=T_\us^\e(X)$. One easily checks that $W$ is a stable algebraic curve and   $T_\us^\e: X\ra W$ is a biholomorphic map.
In particular, we have the following simple lemma.

\begin{lemma}\label{iso}
 If $X_0$ and $X_0'$ are stable analytic curves, and $\us, \us'$ are orthonormal bases for 
$H^0(X_0, mK_{X_0})$ and $H^0(X_0, mK_{X_0'})$ respectively, then the
following conditions are equivalent
\begin{enumerate}
\item $X_0\ \approx \ X_0'$  (i.e. $X_0$ and $X_0'$ are biholomorphic).
\item $[T^\e_{\us'}(X_0')]\ \in \ U(N+1)\cdot [T^\e_\us(X_0)]$
\item $[T^\e_{\us'}(X_0')]\ \in \ SL(N+1,\C)\cdot [T^\e_\us(X_0)]$
\end{enumerate}
Here $[T_\us^\e X_0]\in \hilb$  is the point representing $T_\us^\e X_0\sub\P^N$ in $\hilb$, the Hilbert scheme.
\end{lemma}
\v

\v
Now let $X_i$ be a sequence of stable analytic curves of genus $g$ (e.g Riemann surfaces of genus $g$). Then a basic theorem of Bers \cite{B} (we shall outline the proof below) says there exists a stable analytic curve $X_\i$ (for a 
precise definition see Definition \ref{staban}) such that after passing to a subsequence, $X_i\ra X_\i$.  By this we mean
$X_i^{\rm reg} \ra X_\i^{\rm reg}$ in the pointed Cheeger-Colding topology (see Definition \ref{CG}). Here, for $1\leq i\leq \i$,  $X_i^{\rm reg}\sub X_i$ is the smooth locus. 
This provides the analogue of the left vertical arrow in (\ref{diag}).
 
\begin{theorem}\label{main} Let $X_i$ be a sequence of stable analytic curves of genus $g$. After passing to a subsequence we have 
 $X_i\ra X_\i$ in the Cheeger-Colding sense as above. Then there is a stable algebraic curve $W_\i$ and orthonormal bases $\us_i$ of $H^0(X_i, mK_{X_i})$, such that if $W_i=T_i^\e(X_i)$ then $W_i\ra W_\i$ in the algebraic sense, i.e. as points in the Hilbert scheme. Moreover,  $T_\i|_{X_i^{\rm reg}}$  satisfies 
property (\ref{Tinf}).
\end{theorem}

The idea of using Teichmuller theory  to understand moduli space was advocated by Bers \cite{B, B1, B2, B3}
 in a  of project he initiated, and  which was later completed by Hubbard-Koch \cite{HK}. They define an analytic quotient of ``Augmented Teichmuller Space" whose quotient by the mapping class group is  isomorphic to compactified moduli space as analytic spaces. Our approach is different and is concerned with the imbedding of the universal curve into projective space.
\v
\begin{remark}. As we vary $\e$, the  maps $T_i^\e$ differ by uniformly bounded transformations. We shall
see that if $0<\e_1,\e_2<\e(g)$ then $T^{\e_1}_i=g_i\circ T^{\e_2}_i$ where the change of
basis matrices $g_i\in GL(N+1,\C)$ converge: $g_i\ra g_\i\in GL(N+1,\C)$. In particular,
$\lim_i T^{\e_1}_i(X_i)$ and $\lim_i T^{\e_2}_i(X_i)$ are isomorphic.
\end{remark}
\v
As a corollary of our theorem we shall give a ``metric" proof of the
stable reduction theorem due to  Deligne-Mumford \cite{DM}:
\v
\begin{theorem}\label{1point2} Let $C$ be a smooth curve and  $f:\cX^0\ra C^0$ be a flat family of stable analytic curves over a Zariski open subset $C^0\sub C$. Then there exist a branched cover $\ti C\ra C$ and a flat family $\ti f:\ti\cX\ra \ti C$ of stable analytic curves extending $\cX^0\times_{\ti C}C^0$. Moreover, the extension is unique up to finite base change.
\end{theorem}

 In addition we show that the central
fiber  can be characterized as the Cheeger-Colding limit of the general fibers.
More precisely:
\begin{proposition} Endow $X_t$ with its unique K\"ahler-Einstein metric normalized so that  $\Ric(\o_t)=-\o_t$. Then for every $t\in C^0$ there exist points $p_t^1,....,p_t^\m\in X_t:= f^{-1}(t)$ such that the pointed Cheeger-Colding limits  
$Y_j=\lim_{t\ra 0}(X_t,p_t^j)$ are the connected components of $\ti X_0\backslash\Sigma$ where
$\ti X_0:=\ti f^{-1}(0)$ and $\Sigma\sub \ti X_0$ is the set of nodes of $\ti X_0$. Moreover the limiting
metric on $X_\i$ is its unique K\"ahler-Einstein metric.
\end{proposition}
\begin{remark} A slightly modified proof also gives the log version of stable reduction, i.e for  families $(X_t,D_t)$ where $D_t$ is an effective divisor supported on $n$ points and $K_{X_t}+D_t$ is ample. We indicate which modifications are necessary at the end of section 3.
\end{remark}
\v
\begin{remark} In \cite{S} and \cite{SSW}, Theorems 1.1 and Corollary 2.1 are shown to hold for smooth canonical models of dimension $n>1$. But these papers assume the general version of Theorem 1.2, i.e. of stable reduction. In this paper we do not make these assumptions. In fact, our main purpose here is to prove these algebraic geometry results using analytic methods. 
\end{remark}
\v

We shall first prove Theorem 1.1 under the assumption that the $X_i$ are smooth, and Theorem 1.2 under the assumption that the generic fiber of $f$ smooth. Afterwards we will  treat the general case.

\section{Background}
Let $X$ be a compact connected Hausdorff space, let $r\geq 0$  and $\Sigma=\{z_1,...,z_r\}\sub X$. We say that $X$ is a nodal
analytic curve if $X\backslash \Sigma$ is a disjoint union $Y_1\cup \cdots \cup Y_\m$ of punctured compact Riemann surfaces and if for every $z\in \Sigma$, there is a small open set $z\in U\sub X$ and a continuous function

$$ f: U \ra \{(x,y)\in \C^2\,:\, xy=0\}
$$
with the properties:
\begin{enumerate}
\item $f(z)=(0,0)$
\item $f$ is a homeomorphism onto its image
\item $f|_{U\backslash \{z\}}$ is holomorphic 
\end{enumerate}

\noindent
If $r=0$ then $X$ is a compact Riemann surface.

\v
\begin{definition}\label{staban}
We say that a nodal analytic curve $X$ is a stable analytic curve if each of the $Y_j$ is covered by the Poincar\'e disk.  In other words, each of the $Y_j$ carries a unique hyperbolic metric (i.e. a metric whose curvature is $-1$) with finite volume.
\end{definition}
\v
If $X$ is a stable analytic curve we let $K_X$ be its canonical bundle. Thus the restriction of $K_X$ to $X\backslash\Sigma$ is the usual canonical bundle. Moreover, in the neighborhood of a point $z\in\Sigma$, that is in a neighborhood of
of $\{xy=0\}\sub\C^2$,  a section of $K_X$ consists of a pair of meromorphic differential forms $\eta_1$ and $\eta_2$ defined on $x=0$ and $y=0$ respectively, with the following properties: both are holomorphic away from the origin, both have at worst simple poles at the origin, and $res(\eta_1)+res(\eta_2)=0$. 
\v
We briefly recall the proof of the above characterization of $K_X$ for nodal curves. A nodal singularity is \def\spec{{\rm Spec}} $\spec(B)$  where  $B=\C[X,Y]/(Y^2-X^2)$. Then $\C[X] \ra \C[X,Y]$ is generated by $Y$ which satisfies the monic equation $Y^2-X^2 = 0 $. According the Lipman's characterization of the canonical sheaf \cite{Lip} if $B=C[Y]/(f)$
where $C=\C[X_1,...,X_n]$ and $f$ is a monic polynomial in $Y$ with coefficients in $C$, and if $X=\spec(B)$, then $K_X$ is the
sheaf of holomorphic $(n,n)$ forms on $X_{\reg}$ which can be written as $F\cdot {\pi^*(dx^1\wedge\cdots dx^n)\over f'(Y)}$
where $\pi:X\ra\spec(C)$ and $F$ is a regular function on $X$. In our case, $f(Y)=Y^2-X^2$ so $f'(Y)=2Y$ which means that
$K_X$ is free of rank one, generated by ${dx\over 2y}$ or equivalently ${dx\over y}$. If we consider the map
$\C\ra X$ given by $t\mapsto (t,t)$ then $dx\over y$ pulls back to $dt\over t$. On the other hand, if we consider
$t\mapsto (t,-t)$ then $dx\over y$ pulls back to $-{dt\over t}$.
\v
If $X$ is a compact Riemann surface of genus $g\geq 2$, then $\vol(X)=2g-2$. If $X$ is a stable analytic curve, we say that $X$ has genus $g$ if $\sum_j\vol(Y_j)=2g-2$. Here the volumes are measured with respect to the hyperbolic metric and the $Y_j$ are the irreducible components of $X^{\rm reg}$.
\v
Let $X$ be a stable analytic curve.
The following properties of  $K_X$ are proved in Harris-Morrison \cite{HM}:  \begin{enumerate}
\item $h^0(X,mK_X)=(2m-1)(g-1):=N_m-1$ if $m\geq 2$.

\item $mK_X$ is very ample if $m\geq 3$
\item If $m\geq 3$ the $m$-pluricanonical imbedding of $X$ is a stable algebraic curve in $\P^{N_m}$

\end{enumerate}\v

Next we recall some basic results from Teichmuller theory.
Fix $g>0$ and fix  $S$, a smooth surface of genus $g$.  Teichmuller space $\cT_g$
is the set of equivalence classes of  pairs $(X,f)$ where $X$ is a compact Riemann surface of genus $g$ and $f: S\ra X$ is a diffeomorphism. Two pairs $(X_1,f_1)$ and $(X_2,f_2)$ are equivalent if there is a bi-holomorphic map $h: X_1\ra X_2$ such that 
$f_2^{-1}\circ h\circ f_1: S\ra S$ is in $\Diff_0(S)$, diffeomorphisms isotopic to the identity. The pair $(X,f)$ is called a ``marked Riemann surface". The space $\cT_g$ has a natural topology: A sequence $\tau_n\in\cT_g$ converges to $\tau_\i$ if we can
find representatives $f_n: S\ra X_n$, $1\leq n\leq \i$ such that  the sequence of diffeomorphisms 
$f_\i^{-1}\circ h\circ f_n$ converges to the identity. 
\v

The space $\cT_g$ has a manifold structure given by Fenchel-Nielsen Coordinates whose construction we now recall. Choose a graph $\G$ with the following properties: $\G$ has $g$ vertices, each vertex is connected to three edges (which are not necessarily distinct since  we allow an edge to connect a vertex to itself). For example, if $g=2$, then there are two such graphs: Either $v_1$ and $v_2$ are connected by three edges, or they are connected by one edge, and each connected to itself by one edge.
\v
Fix such a graph $\G$. It has $3g-3$ edges. Fix an ordering $e_1,...,e_{3g-3}$ on the edges. Once we fix $\G$ and we fix an edge ordering, we can define a map 
$(\R_{+}\times\R)^{n} \ra \cT_g$ as follows. Given
$(l_1,\t_1,...,l_n,\t_n)\in \R^{2n}$ we associate to each vertex $v\in\G$ the pair of pants whose geodesic boundary circles have lengths  $(l_i,l_j,l_k)$ where $e_i,e_j,e_k$ are the three edges emanating from $v$. Each of those circles contains two canonically defined points, which are the endpoints of the unique geodesic segment joining it to the other geodesic boundary circles. 
\v
If all the $\t_j=0$, then we join the pants together, using the rules imposed by the graph $\G$, in such a way that canonical points are identified. If some of the $\t_j$ are non-zero, then we rotate an angle of $l_j\t_j$ before joining the boundary curves together. 
\v
Thus we see that $\cT_g$ is a manifold which is covered by a finite number of coordinate charts corresponding to different graphs $\G$ (each 
diffeomorphic to $(\R_{+}\times\R)^{n}$)
If we allow some of the $l_j$ to equal zero, then we can still glue the pants together as above, but this time we get a nodal curve. In this way, 
$(\R_{\geq 0}\times\R)^{n} $ parametrizes all  stable analytic curves.
\v

Teichmuller proved that the  manifold $\cT_g$ has a natural complex structure, and
that there exists a universal curve $\cC_g\ra\cT_g$, which is a map between complex manifolds, such that the fiber above $(X,f)\in \cT_g$ is isomorphic to $X$. Moreover, if
$\cX\ra B$ is any family of marked Riemann surfaces, then there exists
a unique holomorphic map $B\ra \cT_g$ such that $\cX$ is the pullback of 
$\cC_g$. Fenchel-Nielsen coordinates are compatible  with the complex structure, i.e. they  are smooth, but not holomorphic (although they are real-analytic).
\v
\begin{remark}\label{consequence}{One consequence of Teichmuller's theorem is the following. Let
$\cX\ra B$ be a holomorphic family of marked Riemann surfaces and let
$F: B\ra (\R_+\times R)^n$ be the map that sends $t$ to the Fenchel-Nielsen coordinates of $X_t$. Then $F$ is a smooth function. In particular, $X_t\ra X_0$}. This shows that in the stable reduction theorem, if a smooth fill-in exists then it is unique.
\end{remark}

Now let $X$ be a compact Riemann surface. A theorem of Bers \cite{B}, Theorem 15  (a sharp version appears in Parlier \cite{P}, Theorem 1.1) says that for $g\geq 2$ there exists a constant $C(g)$, now known as the Bers constant, with the following property. For every Riemann surface $X$ of genus $g$ there exists a representative $\tau=(X,f)\in\cT_g$ and a graph $\G$ (i.e. a coordinate chart) such that the Fenchel-Nielsen coordinates of $\tau$ are all bounded above by $C(g)$. This is analogous to the fact that  $\P^N$ is covered
by $N+1$ coordinate charts, each biholomorphic to $\C^N$, and that give a point $x\in\P^N$ we can choose a coordinate
chart so that $x\in \C^N$ has the property $|x_j|\leq 1$ for all $j$. In particular, this proves $\P^N$ is sequentially compact.

\v

\v
Bers \cite{B} uses the existence of the Bers constant to show that the space of stable analytic curves is compact with respect to a natural topology (equivalent to the Cheeger-Colding topology). For the convenience of the reader, we recall the short argument.
Let $X_j$ be a sequence of Riemann surfaces. Then after passing to a subsequence, there is a graph $\G$ and representatives $\tau_j=(X_j,f_j)\in\cT_g$  such that the Fenchel-Nielsen coordinates of $\tau_j$ with respect to $\G$ are all bounded above by $C(g)$ (this is due to the fact that there are only finite many allowable graphs). After passing to a further subsequence, we see  $\tau_j\ra\tau_\i\in (\R_{\geq 0}\times\R)^{n} $. If 
$\tau_\i\in (\R_{+}\times\R)^{n}$ then the limit is a smooth Riemann surface. Otherwise, it is a stable analytic curve $X_\i$. Thus

\be\label{union}
X_\i= \cup _{i=1}^\m X^\al,\ \ {\rm and} \ \ X^{\rm reg}_\i=\sqcup_{\al=1}^\m Y^\al
\ee   
where the second union is disjoint, and  $Y^\al = X^\al\backslash F^\al$ where $ X^\al$ is a compact Riemann surface and $F^\al\sub X^\al$ a finite set, consisting of the cusps.
\v

\begin{corollary}\label{cor1} Let $p_\i^\al\in Y^\al$. Then there exist $p_i^1,....,p_i^\m\in X_i$ such that in the pointed Cheeger-Colding topology,
$(Y^\al,p_\i^\al)=\lim_{j\to\i}(X_j,p_j^\al)$. Moreover, for every open set $p_\i^\al\in U_\i^\al\sub\sub Y^\al$ there
exist open sets $p_i^\al\sub U_i^\al\sub X_i$ and diffeomorphisms $f_j^\al: U_\i^\al\ra U_j^\al$ so that $(f_j^\al)^*\o_j^\al \ra \o_\i^\al$ and $(f_j^\al)^*J_j^\al\ra J^\al_\i$ where $\o_j^\al$ and $\o_\i^\al$ are the hyperbolic metrics on $U^\al_j$ and $U_\i^\al$, and 
$J_j^\al$ and $J_\i^\al$ are the complex structures on $U^\al_j$ and $U_\i^\al$
\end{corollary}
\v
\begin{definition}\label{CG} In the notation of Corollary \ref{cor1},
we shall say $\o_j\ra\o_\i$ in the pointed Cheeger-Colding sense and we shall write $X_i\ra X_\i$.
\end{definition}
\v
\noindent
Remark:  Odaka \cite{O2} uses pants decompositions to construct a ``tropical compactification" of moduli space which attaches metrized graphs (of one real dimension) to the boundary of moduli space. These interesting compactifications are compact Hausdorff topological spaces but are no longer algebraic varieties.

\section{Limits of Bergman imbeddings}

\v
Now let $\cX$ be as in the theorem, and let $t_i\in C^0$ with $t_i\ra 0$. Let $X_i=X_{t_i}$  and fix a pants decomposition of $X_i$.
Then Bers' theorem implies that after passing to a subsequence we can find a nodal curve $X_\i$ as above so that
$X_j\ra X_\i$. 
\v

In order to prove the theorem, we must show: 
\begin{enumerate}
\item  $X_\i$ is independent of the choice of subsequence.
\item  After making a finite base change, we can insert $X_\i$ as the central fiber in such a way that the completed family is algebraic.

\end{enumerate}

We begin with (2). Let $X$ be a hyperbolic Riemann surface with finite area (i.e. possibly not compact, but only cusps). The Margulis ``thin-thick decomposition" says that there exists $\e(g)>0$ with the following property. There exists at most $3g-3$ closed geodesics of length less that $\e(g)$. Moreover, for every $\e\leq\e(g)$ the set  

$$ {X\backslash X_\e}\ := \ \{x\in X\,:\, {\rm inj}_x<\e\}
$$
is a finite union of of holomorphic annuli (which are open neighborhoods of short geodesics) if $X$ is compact, and a finite union of annuli as well as punctured disks, which correspond to cusp neighborhoods if $X$ is has singularities. We call these annuli ``Margulis annuli".
Moreover,  $V(\e)$, the volume of $X\backslash X_\e$, has the property $\lim_{\e\ra 0}V(\e)=0$. An elementary proof is given in Proposition 52, Chapter 14 of Donaldson \cite{D}.
\v

Now we define a modified Bergman kernel as follows: For convenience we write $\e=\e(g)$. This is a positive constant, depending only on the genus $g$. Let $X$ be a stable analytic curve. For
$\eta_1,\eta_2\in H^0(X,mK_X)$ let

\be\label{berg1} \<\eta_1,\eta_2\>_\e\ = \ \I_{X_\e} \eta_1\bar\eta_2 h_{KE}^m \o_{KE}
\ee
\noindent
and $\|\eta\|_\e^2=\<\eta,\eta\>_\e$.
If we replace $X_\e$ by $X$, we  get the standard Bergman inner product.

\v

\v

Now fix $m\geq 3$. Choosing orthonormal bases with respect to the inner product (\ref{berg1}) defines imbeddings $T_i^\e: X_i\ra \P^{N_m}$ and $T_\i^\e: X_\i\ra\P^{N_m}$, which we call $\e$-Bergman imbeddings. 
Our goal is to show
\begin{theorem}\label{Tilimit} Let $X_1, X_2, ... $ be a sequence of compact Riemann surfaces of genus $g$. Then there exists a stable analytic curve $X_\i$ such that after passing to a subsequence, $X_i\ra X_\i$ in the Cheeger-Colding topology.  For $1\leq i<\i$, we fix an orthonormal basis $\us_i$ of $H^0(X_i, mK_X)$. Then there exists a choice of orthonormal basis $\us_\i$ for $X_\i$ such that
after passing to a subsequence,
\be\label{limit}  \lim_{i\to\i}T_{i,\us_i}^\e\ = \ T_{\i,\us_\i}^\e
\ee
In other words, if $x_i\in X_i$ and $x_\i\in X_\i$ with $x_i\ra x_\i$, then  $$T_i^\e(x_i)\ra T_\i^\e(x_\i)$$
\end{theorem}\noindent
The proof of Theorem \ref{Tilimit} rests upon the following.
\v
\begin{theorem}\label{ee2} Fix  $g\geq 2$  and $m,\e>0$. Then there exist $C(g,m,\e)$ with
the following property.
$$ \|s\|_\e\ \leq \ \|s\|_{\e/2}\ \leq \ C(g, m,\e)\|s\|_\e
$$
for all Riemann surfaces $X$ of genus $g$ and all $s\in H^0(X,mK_X)$. 
\end{theorem}

To prove the theorem, we need the following adapted version of a
result of Donaldson-Sun. We omit the proof which is very similar to \cite{DS} (actually easier since the only singularities of $X_\i$ are nodes so the pointed limit of the $X_i$ in the Cheeger-Colding topology is smooth).

\begin{proposition}\label{propDS}
Let $X_i\ra X_\i$ be a sequence of Riemann surfaces of genus $g$ converging in the pointed Cheeger-Colding sense to a stable curve $X_\i$. 
Fix   
$\{s_0^\i,...,s_M^\i\}\ \sub\  
H^0(X_\i,mK_{X_\i})$ 
 an $\e$-orthonormal basis   of the bounded sections (i.e. the $L^2(X_\i)$ sections, i.e. the sections which vanish at all nodes). Then there exists an $\e$-orthonormal subset
$$\{s_0^i,...,s_M^i\}\ \sub\  
H^0(X_i,mK_{X_i})$$ 
such that for $0\leq \al\leq M$, we have 
$$s_\al^i\ra s_\al^\i
 $$ 
in $ L^2$
and uniformly on compact subsets of $X_\i^{\rm reg}$. In particular, if $x_i\in X_i^{\rm reg}$
 
\be\label{DS}
\hbox{
$x_i\ra x_\i \iff 
s_\al^i(x_i)\ra s_\al^\i(x_\i)\ \ {\rm for\ all}\ \ 0\leq\al\leq M
\ \ \iff\ \ T^{\n,\e}_i(x_i)\ra T^{\n,\e}_\i(x_\i)$ .
}
\ee
where $T^{\n,\e}_i: X_i^{\rm reg}\hookrightarrow \P^M$ is the map $x_i\mapsto 
(s_0^i,...,s_M^i)(x_i)$ for $1\leq i\leq \i$.

\end{proposition}

\v
{\it Proof of theorem \ref{ee2}.} Let $X_i\ra X_\i$ as in Proposition \ref{propDS}. Choose $(s_0^\i,...,s_M^\i)$ and 
$(t_0^\i,...,t_M^\i)$ which are $\e$ and $\e/2$ orthonormal bases of the subspace of bounded sections in $H^0(X_\i, mK_{X_\i})$ in such a way that
$t_\al^\i=\l_\al^\i s_\al^\i$ for real numbers $0<\l_\al^\i<1$. Choose 
$s_\al^i\ra s_\al^\i$ and $t_\al^i\ra t_\al^\i$ as in Proposition \ref{propDS} in such a way that $t_\al^i=\l_\al^is_\al^i$ with
$0<\l_\al^i<1$. Clearly
\be\label{22} \l_\al^i \ \ra \ \l_\al^\i>0\ \ {\rm for} \ \ 0\leq \al\leq M
\ee
Choose additional sections $s_\al^i$ and $t_\al^i$ for $M+1\leq \al\leq N$
so that $\{s_0^i,..., s_N^i\}$ and $\{t_0^i,..., t_N^i\}$ are $\e$ and $\e/2$ bases of $H^0(X_i,mK_{X_i})$ and $t_\al^i=\l_\al^is_\al^i$ with
$0<\l_\al^i<1$ for $0\leq\al\leq N$.
\v
Now assume the theorem is false. Then there exists $X_i\ra X_\i$ as above such that $\l_\al^i\ra 0$ for some $\al$.
 We must have $\al\geq M+1$ by (\ref{22}).
Choose $M+1\leq A<N$ such that $\l_\al^i\ra 0$ if and only if $A\leq \al\leq N$.
\v\noindent
Since $\|s_\al^i\|_{L^2(X_\e)}=1$ we may choose $s_\al^\i(\e)\in H^0(X_i^\e, K_{X_\i}|_{X_\i^\e})$ such that
\be s_\al^i|_{X^\e_i} \ra s_\al^\i(\e)\ \ {\rm for}\ \ M+1\leq \al\leq N
\hbox{\ \ uniformly on compact subsets of $X_\e$}
\ee
 
Let $T_i^\e: X_i\ra W_i^\e\sub\P^N$ be the Kodaira map given by the
sections $s_0^i,...,s_N^i$ and let {$W_\i^\e=\lim_{i\to\i}W_i^\e$. }
Let
$$ T_\i^\e: X_\i^\e \hookrightarrow W_\i^\e \ \ {\rm and} \ \ 
T_\i^{\n,\e}: X^{\rm reg}_\i \hookrightarrow\P^M
$$
be the Kodaira maps given by $(s_0^\i,...,s_M^\i, s_{M+1}^\i(\e),...,s_N^\i(\e))$ and $(s_0^\i,...,s_M^\i)$.  Thus
\be\label{inc} \pi\circ T_\i^\e\ = \ T_\i^{\n,\e}|_{X_\i^{\e}}
\ee
where
$\pi:\P^N_M:=\P^N\backslash\{z_0=\cdots = z_M=0\}\ra \P^M$ is
defined by $(z_0,...,z_N)\mapsto (z_0,...,z_M)$. Moreover
$$ \pi(W^\e_\i\cap\P^N_M)\ \sub \ T_\i^{\n,\e}(X_\i^{\rm reg})
$$

\v
Now the definition of $A$ implies
$$  T_\i^{\e/2}(X^{\e}_{\i})\ \sub\ Z^{\e/2}_\i\ := \ \{z\in  W^{\e/2}_\i\,:\, z_{A}=z_{A+1}=\cdots=z_N=0\}
$$
Thus (\ref{inc})  implies
$$ T_\i^{\n,\e/2}(X_\i^{\rm reg})\ \supset\ 
\pi(Z^{\e/2}_\i\cap \P^N_M)\supset \pi(T_\i^{\e/2}(X^{\e}_{\i}))\ = \ 
T_\i^{\n,\e/2}(X_\i^\e)
$$
Since the second set is constructible,
$$ \pi(Z^{\e/2}_\i\cap \P^N_M) \ = \  T_\i^{\n,\e/2}(X_\i^{\rm reg}\backslash \Sigma_\e)
$$
where $\Sigma_\e\sub X_\i^{\rm reg}\backslash X_\i^\e$ is a finite set. Moreover, $\Sigma_\e$ is monotone in $\e$.
\v
Let $x_\i\in X_\i^{\rm reg}\backslash\Sigma_\e$. Then $T_\i^{\n,\e/2}(x_\i) = \pi(w_\i)$
for some $w_\i\in Z_\i^{\e/2}\cap\P^N_M$. Choose $w_i\in W_i^{\e/2}$ such
that $w_i\ra w_\i$ and choose $x_i\in X_i$ such that $T_i^{\e/2}(x_i)=w_i$.
Then (\ref{DS}) implies

$$ T^{\e/2}_i(x_i)\ra w_\i\ \Lra\ \pi(T^{\e/2}_i(x_i))\ra \pi(w_\i)\ \Lra \ T_i^{\n,\e/2}(x_i)\ \ra \ T_\i^{\n,\e/2}(x_\i) \ \Lra \ x_i\ra x_\i
$$

Thus we see that if $x_\i\in X_\i^{\rm reg}\backslash\Sigma_\e$ there exists
$x_i\ra x_\i$ such that 

$$ \lim_{i\to\i} T_i^{\e/2}(x_i) \ \in \ Z_\i^{\e/2}
$$
\begin{lemma}\label{Sigma} If $\e_1<\e_2$ then $\Sigma_{\e_1}\sub \Sigma_{\e_2}\cap X_\i\backslash X_\i^{\e/2}$.
\end{lemma}
{\it Proof.} Assume $x_\i\notin \Sigma_{\e_2}$. Then there exists $x_i\ra x_\i$ such that
$ \lim_{i\to\i}s^i_\al(x_i) =0
$
for all $A\leq \al\leq N$. Then
$ \lim_{i\to\i}t^i_\al(x_i) =\lim_{i\to\i}\l_\al^is^i_\al(x_i) =0
$
since $0<\l_\al^i<1$. \qed
\v
We may there assume, after possibly decreasing $\e$, that $\Sigma_\e=\emptyset$.  This means that for every $x_\i\in X_\i^{\rm reg}$
there exists $x_i\ra x_\i$ such that
$ \lim_{i\to\i}s^i_\al(x_i) =0
$
for all $A\leq \al\leq N$.
\v
Let $x_\i\in X_\i^{\rm reg}$. We say that $x_\i$ is an \emph{$\e$-good} point if
for every $x_i\ra x_\i$, $ \lim_{i\to\i}s^i_\al(x_i) =0
$
for all $A\leq \al\leq N$. The set of $\e$-bad points is finite (otherwise
$W^\e_\i$ would have infinitely many components {by the intermediate value theorem}). Also, every point in $X_\i^{2\e}$ is $\e$-good.

\v

\begin{lemma}\label{good} If $x_\i$ is $\e_2$-good then it is 
$\e_1$-good for all $\e_1<\e_2$.
\end{lemma}
{\it Proof.} Assume $x_\i$ is $\e_2$ good, let $x_i\ra x_\i$ and let $\e_1<\e_2$. Then
$ \lim_{i\to\i}t^i_\al(x_i) =\lim_{i\to\i}\l_\al^is^i_\al(x_i) =0
$
since $0<\l_\al^i<1$. \qed
\v
Lemma \ref{good} implies that by decreasing $\e$ if necessary, that all points $x_\i\in X_\i^{\rm reg}$ are $\e$-good. But $\|s_A^i\|_\e=1$ so there exists $x_i\in X_i^\e$ such that $|s_A^i(x_i)|=1$. After passing to a subsequence, $x_i\ra x_\i\in X_\i^\e$ but
$\lim_i|s_A^i(x_i)|=1\not=0$ a contradiction.

\v

We conclude that if $\eta_j\in H^0(X_j, mK_{X_j})$ is a sequence
such that the norms $\|\eta_j\|^2_\e=\<\eta_j,\eta_j\>_\e =1$, then after passing to a subsequence, we have $(f^\al_j)^*\eta_j \ra  \eta_\i$ for some $\eta_\i\in H^0(X_\i^{reg}, mK_{X_\i}|_{X_\i^{\rm reg}})$ with $\|\eta\|_\e=1$. Here the $f_j^\al:U^\al_j\ra U^\al$ are as in the statement of Corollary 1 and this is true for all $U^\al$ and all $\al$.  Moreover, an
orthornormal basis of $H^0(X_j, mK_{X_j})$, which is a vector space of dimension $(2m-1)(g-1)$,  will converge to an orthonormal set of $(2m-1)(g-1)$ elements in 
$H^0(X_\i^{reg},mK_{X_\i})$. The main problem is to now show that these $(2m-1)(g-1)$ elements extend to elements of $H^0(X_\i,mK_{X_\i})$. If they extend, then they automatically form a basis since 
$H^0(X_\i,mK_{X_\i})$ has dimension $(2m-1)(g-1)$ and this would prove
Theorem \ref{Tilimit}.
\v

\v
To proceed, we make use of the discussion of the Margulis collar  in  section 14.4.1 of \cite{D}.  Let $\l>0$ be the length of $C$ a collapsing geodesic in $X_j$ which forms a node in the limit in $X_\i$. We fix $j$ and we write $X=X_j$. Let

\v
$$ A_\l\ = \ \{z\in\C: 1\leq |z|\leq e^{2\pi\l},\, \l\leq \arg(z)\leq \pi-\l\,\}/\sim
$$
where the equivalence relation identifies the circles $|z|=1$ and $|z|=e^{2\pi\l}$.
Then \cite{D} shows $A$ injects holomorphically into $X$ in such a way that $1\leq y\leq e^{2\pi\l}$
maps to $C$. The point is that the segment $1\leq y\leq e^{2\pi\l}$ is very short - it has size $\l$. But the segments $A\cap \{\arg(z)=\l\}$ and $A\cap \{\arg(z)=\pi-\l\}$ have size 1. So for $\l$ small, $A$ is a topologically a cylinder, but metrically very long and narrow in the middle but not narrow at the ends.  In other words, the middle of $A$ is in the thin part, but the boundary curves are in the thick part.

\v


 The transformation

$$ \tau\ = \ \exp\big(i{\ln z\over \l}\big)
$$
maps $A_\l$ to the annulus 

$$ A'_\l \ = \ \{\exp(-(\pi-\l)/\l)\leq |\tau|\leq \exp(-1)\,\}
$$

To summarize: We are given a sequence $X_j$, and a geodesic $C_j$ in $X_j$ that collapses to a node $\nu$ in $Y^\al$ for some $\al$. We are also given a sequence of orthonormal bases
$\{{\eta_{j,1}},..., \eta_{j,N}\}$ of 
$H^0(X_j,kK_{X_j})$ where $N=(2k-1)(g-1)$ and $\eta_{j,\m}\ra \eta_{\i,\m}$. Here $\eta_{\i,\m}$ is a section of $kK_{X\i}$ on $X_\i^{\rm reg}$. Fix $\m$ and write $\eta_j=\eta_{j,\m}$ and $\eta_\i=\eta_{\i,\m}$. We need to show that $\eta_\i$ extends to all of $X_\i$.
\v
We may view $\eta_j$ as a $k$ form on $A_{\l_j}$ or on $A'_{\l_j}$ and $\eta_\i$ as a $k$ form on the punctured disk $A'_0$. Write
$\eta_j = f_j(z)dz^k=h_j(\tau)d\tau^k$ and $\eta_\i=h_\i(\tau)d\tau^k$. The discussion in \cite{D} shows that if we fix a relatively compact open subset  $U\sub A'_0$, then $h_j \ra h_\i$ uniformly on $U$ (this makes sense since $U\sub A'_{\l_j}$ for $j$ sufficiently large).

\v
Since $\|\eta_j\|_{L^2}=1$  we have uniform sup norm bounds on the thick part of $X_j$. Thus

\be\label{thick} \|\eta_j\|_{L^\i((X_i)_\e}\ \leq \ C(\e)
\ee

We want to use (\ref{thick}) to get a bound on the thin part. In $z$ coordinates, (\ref{thick}) implies 

\be\label{sup}  |\eta|_{\o}\ = \ |\Im(z)|^k\cdot |f(z)|\ \leq \ C(\e) \ \ \hbox{if $\arg(z)=\l$ or $\arg(z)=2\pi-\l$}
\ee
since the boundary curves $\arg(z)=\l$ and $\arg(z)=2\pi-\l$ are in the thick part.
Here we write $\eta$ for $\eta_j$ and $f$ for $f_j$ to lighten the notation.
\v
Now

\be\label{23}\Im(z)\ = \ -\exp(\l\arg\tau)(\sin(\l\ln|\tau|)\ee

if we write $f(z)=g(\tau)$, then (\ref{sup}) implies

\be |g(\tau)|\leq {C(\e)\over\l^k}\ \ \hbox{for $\tau\in\pl A'$}
\ee

Since $f(z)dz^k=h(\tau)d\tau^k=g(\tau)({dz\over d\tau})^k\,d\tau^k$ and
${dz\over d\tau}= {z\l\over i\tau}$
we see  for $\l$ small

$$ |h_j(\tau)|\ \leq {1\over \l^k}\bigg|{dz\over d\tau}\bigg|^k\ = \ 
{1\over\l^k}{|z|^k\l^k\over|\tau|^k}\ \leq \ {2\over |\tau|^k}
$$
where the last inequality follows from the fact $1\leq |z|\leq 2$.
Writing 
$$ u_j(\tau)\ :=\ h_j(\tau)\tau^k
$$
Thus we see $|u_j(\tau)|\leq 2$ for $\tau\in\pl A'$. The maximum principle now implies that $ |u_j(\tau)|\leq 2$ for $\tau\in A'$. Since this is true for all
$X_i$, we see that any limit $u_\i$ must satisfy the same inequality in the limit of the annuli, which is a punctured disk: $|h_\i(\tau)|\cdot |\tau|^k\leq C$. This shows $h_\i$ has at most a pole of order $k$. 

Moreover  $u(0)$ is
the residue \be\label{res} u(0)\ = \ \lim_{j\to\i}{1\over 2\pi \sq}\I_{|\tau|=r} u_j(\tau)\,{d\tau\over\tau}
\ee
Here $0<r\leq \exp(-1)$ is any fixed number (independent of $j$).

\v
To summarize, we have now seen that a collar degenerates to a union of two punctured disks and so the limit of
the $\eta_j$ is a pair of $k$ forms, $\eta_\i = u_\i(\tau)\left({d\tau\over \tau}\right)^k$ and 
$\ti\eta_\i = \ti u_\i(\tau')\left({d\tau'\over \tau'}\right)^k$ where $u$ and $\ti u$ are holomorphic
in a neighborhood of the origin in $\C$. There is one final condition that we need to check in order
to verify that the limit is in $H^0(X_\i,kK_{X_\i})$:
 Let $R=\exp(-1)$, $r=\exp(-\pi/2\l_j)$ and $\e= \exp(-\pi/\l_j)$ (so $\e_{\l_j}/r = r$).  We must show $\ti u(0)=(-1)^ku(0)$.
\begin{center}
\begin{figure}[h!]
\begin{tikzpicture}
\coordinate (O) at (0,0);

%

\begin{scope}[xshift=6cm]
\coordinate (O) at (0,0);
\draw[fill=blue!25] (O) circle (2.8);
\draw[fill=red!25] (O) circle (1.2);
\draw[fill=white!25] (O) circle (0.4);

\node[text width=2cm] at (.1,-.65) {\tiny\black ${\e/R}\leq |\tau|\leq r$};
\node[text width=2cm] at (.35, 3.05) {\tiny\black $|\tau|= R$};
\node[text width=2cm] at (.2, -2) {\tiny\black $r\leq |\tau|\leq R$};

\end{scope}
\end{tikzpicture}
\caption{$A_\l'$}\label{X0}
\end{figure}
\end{center}

\v
\noindent
(Here the inner white disk is $|\tau|\leq \e/R$).
To check this, let $\ti\tau= {\e_j\over\tau}$. Then Figure \ref{X0} remains the same, with $\tau$ replaced by $\ti\tau$ and 
$$ f(z)dz^k\ = \ u_j(\tau)\left({d\tau\over\tau}\right)^k\ = \ u_j(\e_j/\ti\tau)(-1)^k\left({d\ti\tau\over\ti\tau}\right)^k
\ := \ \ti u_j(\ti\tau)\left({d\ti\tau\over\ti\tau}\right)^k
$$

Now we see

$$
\I_{|\tau|=r} u_j(\tau)\,{d\tau\over\tau}\ = \  
(-1)\I_{|\ti\tau|=\e_j/r} u_j({\e_j\over\ti\tau})\,(-1){d\ti\tau\over\ti\tau}\ = \ 
(-1)^k\I_{|\ti\tau|=r} \ti u_j(\ti\tau)\,{d\ti\tau\over\ti\tau}
$$
In the second integral, the factor of $(-1)$ outside the integral is due to the fact that the orientation of the 
circle has been reversed and the $(-1)$ inside the integral comes from the change of variables. The second identity
is a result of the fact that $u(\ti\tau)$ is holomorphic on the annulus $\{\ti\tau\in\C:\e_j/r < \ti\tau< r\}$.
Taking limits as $j\to\i$ we obtain $\ti u(0)=(-1)^ku(0)$. This establishes Theorem \ref{Tilimit} when the $X_i$ are smooth.

\v
Now assume the $X_i$ are stable analytic curves, but not necessarily smooth. The Fenchel-Nielsen coordinates of $X_i$ determine a point  
$[X_i]\in (\R_{\geq 0}\times\R)^{n} $. The simple observation we need is that $(\R_{> 0}\times\R)^{n}\sub (\R_{\geq 0}\times\R)^{n} $ is dense
so we may choose a smooth Riemann surface $\ti X_i$ such that $[X_i]
\in (\R_{\geq 0}\times\R)^{n} $ is $\e_i$  close to $[X_i]$ where $\e_i\ra 0$ (i.e. $X_i$ is smoothable). Now Corollary \ref{cor1} implies that after passing to a subsequence,
$\ti X_i\ra X_\i$ in the pointed Cheeger-Colding topology. We conclude
that $X_i\ra X_\i$ as well. Moreover, one easily sees that $T_i^\e$ and 
$\ti T_i^\e$ have the same limit. This proves
(\ref{limit}) and
completes the proof of Theorem \ref{Tilimit}
 \qed
\v
\begin{remark} The proof of the log version Theorem \ref{Tilimit} is almost the same. The only observation  we need is  the following. If $X$ is a compact Riemann surface and $D=p_1+\cdots+p_n$ is a divisor supported on $n$ points such that $K_X+D$ is ample, then $X\backslash D$ has a unique metric $\o$ such that $\Ric(\o)=-\o$ and $\o$ has cusp singularities at the points $p_j$. Morover, just as in the case $n=0$, $X$ has a pants decomposition. The only difference is that we allow some of the length
parameters to vanish, but this doesn't affect the arguments. In particular, we can use the Fenchel-Nielson coordinates to find a limit of the $(X_j,D_j)$ (after passing to a subsequence) and the $T_j^\e$ are defined exactly as before.
\end{remark}
\v
Now suppose $X_i$ is a sequence of compact Riemann surfaces of genus $g$ converging analytically to a nodal curve $X_\i$
and let $\eta_i$ be a \K metric on $X_i$ is the same class as the \KE metric $\o_i$. We have seen that $\o_i\ra\o_\i$,
the \KE metric on~$X_\i$, in the pointed Cheeger-Colding sense. Let 
$\ti\o_\i$ 
be a \K metric on $X_\i^{\reg}$ and assume 
$\ti\o_i\ra \ti\o_\i$ in the pointed Cheeger-Colding sense. Let $T_i(\ti\o_i): X_i\ra \P^N$ be the embedding defined by an orthonormal basis of $H^0(X_i,3K_{X_i})$ using the metric $\ti\o_i$ on the thick part of $X_i$
and define $T_\i(\ti\o_\i):X_\i\ra \P^N$ similarly. Thus the $T_i$ and $T_\i$ of Theorem 2 can be written as
$T_i(\o_i)$ and $T_\i(\o_\i)$ and in this notation, Theorem 2 says 
$T_i(\o_i)\ra T_\i(\o_\i)$
\begin{corollary} After passing to a subsequence
$$ T_i(\ti\o_i)\ \ra \ T_\i(\ti\o_\i)
$$
\end{corollary}

{\it Proof.} Since $\ti\o_\i$ and $\o_\i$ are equivalent on the thick part of $X_\i$, we see that
$$T_i(\ti\o_i)=\g_i\circ T_i(\o_i)$$
where $\g_i\in GL(N+1,\C)$ has uniformly bounded entries as does 
$\g_i^{-1}$. Thus
after passing to a subsequence, $\g_i\ra\g_\i\in GL(N+1,\C)$ and 
$$\lim_{i\to\i} T_i(\ti\o_i)\ = \ \lim_{i\to\i}\g_i\circ T_i(\o_\i)\ = \ \g_\i\circ T_\i(\o_\i)\ = \  T_\i(\ti\o_\i)
$$
\noindent
Remark: The proof shows we only need to assume $\ti\o_i\ra \ti\o_\i$ on the thick part of $X_\i$.

\section{Existence of stable fill-in} \setcounter{equation}{0}
\v
\noindent
{\it Proof.}
Let $f:\cX^0\ra C^0=C\backslash\{p_1,...,p_m\}$ be a flat family of stable analytic curves of genus $g\geq 2$. We shall assume the fibers are smooth since except for some additional noation, the general case is proved in exactly the same way.  We first observe that we can find some completion (not necessarily nodal) 
$\cY\ra C$ of the family $\cX^0\ra C^0$. To see
this let $\O_{\cX^0/C^0}$ be the sheaf of relative differential forms (i.e. the relative canonical line bundle when $\cX^0$ is smooth). Then the Hodge
bundle $f_* K_{\cX^0/C^0}$ is a vector bundle over $C^0$ of rank $3g-3$ (see page 694 of Vakil [V])
and $f_*K_{\cX^0/C^0}^{\otimes m}$ is a vector bundle $\cE_m^0$ of rank $N_m-1:=(2m-1)(g-1)$ for $m\geq 2$.
Choose $\cE_m\ra C$ an extension of the vector bundle  $\cE_m^0\ra C^0$ to the curve $C$.
\v
For example, let $U\sub C^0$ be any affine open subset over which $\cE^0_m$ is trivial and let $s_0,...,s_{N_m}$ be a fixed $\cO(U)$ basis. Then if $p_j\in V\sub C^0$ is an affine open set such that $V\backslash \{p_j\}\sub U$, then define $\cE(V)$ to be the $\cO(V)$ submodule of 
$\cE^0(V\backslash\{p_j\})$ spanned by
the $s_\al$. 
\v
 Once $\cE$ is fixed, we choose $m\geq 3$ and let $\cX^0\hookrightarrow \P(\cE^0)\sub\P(\cE)$ be the canonical imbedding. Then we define 
 \be\label{fake} \cY\sub \P(\cE)
 \ee
  to be the flat limit of 
$\cX^0\ra C^0$ inside $\P(\cE)\ra C$.

\v
Now we prove Theorem 1. To lighten the notation, we shall assume $m=1$ and write $C^0=C\backslash\{0\}$ where $0:= p_1$.
Suppose $t_i\in C^0$ with $t_i\ra 0$ and such that we have analytic convergence  $X_{t_i}\ra X_\i$ where $X_\i$ is an stable analytic curve.
 We wish to show that there exists a smooth curve $\ti C$ and a finite cover $\m:\ti C\ra C$ with the following property. If we let $\Sigma=\m^{-1}(0)$ (a finite set) there exists a unique completion
$\ti f:\ti\cX\ra \ti C$ of $\m^*\cX^0\ra \ti C\backslash\Sigma$ with $X_\i=p^{-1}(\ti 0)$ for all $\ti 0\in \Sigma$.
\v
Define
$$ Z^0\ = \ \{(t,z)\in C^0\times \Hilb(\P^{N_m})\,:\ z\in \cT_t\,\}
$$
where $\cT_t$ is the set of all Hilbert points $[T(X_t)]$. Here $T:X_t\ra\P^{N_m}$ ranges over the set of all Bergman imbeddings. In particular, $\cT_t\sub \Hilb(\P^{N_m})$ is a single  $G=SL(N_m+1)$  orbit.  

\v
We claim that $Z^0\sub C^0\times\Hilb(\P^{N_m})$ is a constructible subset. To see
this, let $U\sub C^0$ be an affine open subset
and let $\si_0,...,\si_{N_m}$ be a fixed $\cO(U)$ basis of $\cE_m(U)$. This basis defines an imbedding

\be\label{vb} S:\pi^{-1}(U)\ \ra \ U\times\P^{N_m}
\ee
given by $x\mapsto (\pi(x), \si_0(x),...,\si_{N_m}(x))$. Define $H:U\ra \Hilb(\P^{N_m})$ by $H(t)= \Hilb(S(X_t))$ and define the map
$$ f_U:G\times U\ \ra \ U\times\Hilb(\P^{N_m})\ \ 
\hbox{given by $(g,t)\mapsto (t, g\cdot H(t))$}
$$
Then $f_U$ is an algebraic map so its image is constructible. This shows $Z^0|_{U}$ is
constructible for every affine subset $U\sub C^0$ and hence $Z^0$ is constructible.
\v
Now we 
fix $0<\e<\e(g)$ and
let $W_j=T_j(X_{t_j})$ where $T_j$ is the $\e$-Bergman imbedding. Then (\ref{limit}) implies  $T_j(X_j)=W_j\ra T_\i(X_\i)= W_\i$, a stable algebraic curve in $\P^{N_m}$.
Let $Z\ra C$ be the closure of $Z^0$ in $C\times \Hilb(\P^{N_m})
\sub C\times\P^M$. Here $\P^M\supset \Hilb(\P^{N_m})$ is chosen so that there is a $G$ action on $\P^M$ which restricts to the $G$ action on  $\Hilb(\P^{N_m})$. Then $Z$ is a subvariety of $C\times \Hilb(\P^{N_m})$ whose dimension we denote by $d$. Let $Z_t$ the fiber of $Z$ above $t\in C$.
 Then $[Y_\i]\in Z_0$. 
 \v
To construct $\ti C$ we use the Luna Slice Theorem: 
There exists $W\sub\C^{M+1}$ a $G_{[Y_0]}$ invariant subspace such that 
$[Y_\i]\in \P(W)\sub\P^M$
and such that the map
$$ \P(W)\times \Lie(G)\ \ra \ \P^{M}\ \ \hbox{
given by\ $(x,\xi)\mapsto \exp(\xi)x$}
$$
is a diffeomorphism of some small neighborhood $U_W\times V\sub \P(W)\times\Lie(G)$ onto an open set $\O\sub\P^M$,  with $U_W\sub\P(W)$ invariant under the finite group $G_{[Y_0]}$. After shrinking $U_W$ if necessary, the
intersection of a $G$ orbit with $U_W\backslash [Y_0]$ is a finite set of  order $m_1|m$ where
$m=|G_{[Y_0]}|$. In other words, the quotient $G_{[Y_0]}\backslash U_W$ parametrizes
the $G$-orbits in $\P^M$ that intersect $U_W$.
\v
Note that $\O$ contains $(t_i, [Y_i])$ for infinitely many $i$ so $(C\times \P(W))\cap Z$ is a projective variety $C_1$ of dimension at least one. Moreover, if we let $C_2$ be the union of the components of $C_1$ containing $\{0\}\times [Y_\i]$, then $C_2$ contains infinitely many of 
$(t_i, [Y_i])$ so the image of $C_2\ra C$ contains infinitely many $t_i$ and thus $C_2\ra C$ is surjective. On the other hand, $C_2\ra C$ is finite of degree $m_1$ (this follows from the construction of $U(W)$).
\v
Let $\ti C\sub C_1$ be an irreducible component of $C_1$ containing
$(t_i, [Y_i])$ for infinitely many $i$. 
Let $H\sub G_{[Y_\i]}$ be the set of all $\sigma\in 
G_{[Y_\i]}$ such that $\sigma(\ti C)=\ti C$. Then $H$ has order $d$ for
some $d|m_1$ and $\ti C\ra C$ is finite of degree $d$.
\v
Finally, we have $\ti C\sub Z\sub C\times \Hilb(\P^{N_m})$. This gives us
a map $\ti C\ra \Hilb(\P^{N_m})$. If we pull back the universal family
we get a flat family $\ti \cX\ra \ti C$ which extends
$\cX^0\times_{\ti C}C^0$. This completes the proof of Theorem 1.2.\qed

\v

\section{Uniqueness of the stable fill-in} \setcounter{equation}{0}
 Let $\pi:X^*\ra\D^*\subset \D$ be an algebraic     family of smooth  curves genus $g$. We claim that
there exists a  unique stable analytic curve $X_0$ such that $X_t\ra X_0$ in the Cheeger-Colding sense
as $t\ra 0$.  This will establish the uniqueness statement of Theorem \ref{1point2}, and since existence was demonstrated in the previous section, it completes the proof.

\v
Let $S: \cX^*\ra \D^*\times\P^{N_m}$ as in (\ref{vb}). For each $t\in\D^*$, the set $\u\si_t=(\si_0(t),...,\si_{N_m}(t))$ is a basis of $H^0(X_t, mK_{X_t})$.
Let $\us_t=(s_0(t),...,s_{N_m}(t) )$ be the orthonormal basis
of $H^0(X_t, mK_{X_t})$ obtained by applying the Gram-Schmidt process
to  the basis $\u\si_t$ and let $T^\e_t: X_t\ra\P^N$ be the map $T^\e_t=T^\e_{\us_t}$.
Here $0<\e<\e(g)$ is fixed once and for all. Remark \ref{consequence} implies that $t\mapsto [T_t^\e(X_t)]$ defines a continuous function
$\D^*\ra \hilb$. Let

$$ z: \D^*\times SL(N+1,\C)\ \ra \ \D^*\times\hilb
$$
and 
$$f: \D^*\ra \D^*\times \hilb$$
be the maps
$$ z(t,g)\ = \ (t, g\cdot [T_t(X_t)])\ \ 
\hbox{and $f(t)\ = \ z(t, [T_t(X_t)])$.}
$$

\v

Let $F=\overline{\Im(f)}\sub \D\times \hilb$ and $Z=\overline{\Im z}\sub\D\times\hilb$. Let $\pi_F:F\ra \D$ and $\pi_Z:Z\ra\D$
be the projection maps and $F_0=\pi_F^{-1}(0)$, $Z_0=\pi_Z^{-1}(0)$.
Observe that $F_0\sub\hilb$ is closed and connected (this easily follows from the fact that $\D^*$ is connected and $\hilb$ is compact and connected). Moreover, Theorem \ref{Tilimit} implies that every element of $F_0$ is of the form $T^\e_\us(X_0)$ for some stable analytic curve $X_0$ and some basis $\us$.
\v
\noindent
Claim:  $F_0$ is contained in the $U(N+1)$ orbit of $[X_0]$.
\v
Assume the claim for the moment, and let's show that it implies uniqueness.
Suppose there exist subsequences $t_i, t_i'\in \D^*$ such that 
$X_{t_i} \ra X_0$ and $X_{t_i'}\ra X_0'$. We must show that $X_0\approx X_0'$, i.e. $X_0$ and $X_0'$ are isomorphic stable analytic curves. Theorem \ref{Tilimit} implies there are bases $\us$ and $\us'$ such that $[T^\e_\us(X_0)], [T^\e_{\us'}(X_0')]\in F_0$ so 
$ T^\e_{\us'u}(X_0')\in U(N+1)\cdot T^\e_\us(X_0)$. Now Lemma \ref {iso} implies  $X_0\approx X_0'$. This gives uniqueness.

\v
\noindent
The set  $U=SL(N+1,\C)\cdot [T^\e_\us(X_0)]\sub Z_0$ is open since $\dim Z_0=\dim SL(N+1,\C)$ and the stabilizer of $[T^\e_\us(X_0]$ is finite.
Lemma \ref{iso} implies 
\be\label{U} F_0\cap U\ \sub \ U(N+1)[T_{\us}^\e(X_0)] \ \sub \ U
\ee
Now $U(N+1)[T_{\us}^\e(X_0)]$ is compact and $F_0$ is connected, so
$F_0\cap U = F_0$. Thus the claim follows from (\ref{U}). \qed
\v

\end{document}